\documentclass[12pt, leqno]{article}
\usepackage{amsmath}
\usepackage{amsfonts}
\usepackage{amssymb}
\usepackage{amssymb,esint}
\usepackage{hyperref}
\usepackage{graphicx}
\usepackage{color}
\usepackage{ifpdf}     
\usepackage{pdfsync}   

\newtheorem{thm}{Theorem}[section]

\newtheorem{lem}[thm]{Lemma}
\newtheorem{prop}[thm]{Proposition}

\numberwithin{equation}{section}

\newcommand{\RR}{\mathbb{R}}

\newcommand{\ve}{\varepsilon}

\newcommand{\ren}{{\RR^N}}
\newcommand{\s}{{\rho}}

\def\qed{\,\unskip\kern 6pt \penalty 500
\raise -2pt\hbox{\vrule \vbox to8pt{\hrule width 6pt
\vfill\hrule}\vrule}\par}
\definecolor{darkblue}{rgb}{0.05, .05, .65}
\definecolor{darkgreen}{rgb}{0.1, .65, .1}
\definecolor{darkred}{rgb}{0.8,0,0}

\begin{document}
\title{\textbf{The Mesa Problem for the Fractional \\[8pt] Porous Medium Equation}\\[7mm]}

\author{\Large Juan Luis V\'azquez\\}

\date{}

\maketitle

\begin{abstract}
We investigate the behaviour of the solutions $u_m(x,t)$ of the  fractional porous medium equation
$$
u_t+(-\Delta)^s (u^m)=0, \quad x\in \ren, \ t>0.
$$
with initial data $u(x,0)\ge 0$, $x\in \ren$, in the limit as $m\to\infty$ with fixed $s\in (0,1)$. We first identify the limit of the Barenblatt solutions as the solution of a fractional obstacle problem, and we observe that, contrary to the case $s=1$, the limit is not compactly supported but exhibits a typical fractional tail with power-like decay. In other words, we do not get a plain mesa in the limit, but a mesa with tails. We then  study the limit for a class of nonnegative initial data and derive counterexamples to expected propagation and comparison properties based on symmetrization.

\end{abstract}


\



\section{The mesa problem as limit of the FPME when $m\to\infty$}
\setcounter{equation}{0}

This paper deals with the limit $m\to\infty$ of the solutions of the family of fractional porous medium equations
(FPME) of the form
\begin{equation}\label{eq.ms}
u_t+(-\Delta)^s (u^m)=0, \quad x\in \ren, \ t>0.
\end{equation}
Here $(-\Delta)^s$ is the fractional Laplacian operator with $0<s<1$. We will take initial data
\begin{equation}\label{in.data}
u(x,0)=u_0(x),
\end{equation}
where $u_0$ is a nonnegative and integrable real function defined in $\ren$, or a Dirac delta. Under the former  assumptions this Cauchy Problem  produces a unique solution $u_m(x,t)$ for every finite $m>1$ and every $0<s<1$, cf. \cite{pqrv, pqrv2} for the basic theory and the survey paper \cite{VazMilan} for recent progress in the qualitative analysis. As $s\to 1$ with $m>1$ fixed we obtain the limit the standard porous medium equation (PME), \ $u_t=\Delta u^m$, whose mathematical theory  and applications are described in  \cite{vazquezPME}. In the case of a Dirac delta as initial data, the solutions are called  fundamental solutions or Barenblatt solutions; in the fractional case  the Barenblatt solutions have been constructed in \cite{vazBaren} and their uniqueness, properties and role in the asymptotic behaviour are explained, see Section \ref{sec.barfrac} below.

The study of the differences in the behaviour of  diffusion equations involving the standard Laplacian, i.\,e., involving a local operator, and the fractional variants, which involve nonlocal diffusion operators,  is a matter of much current attention. The comparison between the PME and the FPME is a convenient setting for such analysis. We tackle here the limit $m\to\infty$ in both models.

The dramatic effect of the limit $m\to\infty$ is well known in the PME case, i.\,e., equation \eqref{eq.ms} with $s=1$, as described in \cite{CaffFr87, EHKO86, FrHo87, FrHu88, Sacks1989, BBH89, BI2003, BI2004, Igb2002} and other related literature.  It is proved that given a nonnegative initial datum $u_0(x)\in L^1(\ren)$ there is a unique weak solution $u_m(x,t)$ of the PME defined in $Q=\ren\times(0,\infty)$ that is uniformly bounded and $C^\alpha$ continuous for $t\ge s>0$. If $u_0$ is compactly supported, so is $u_m(\cdot,t)$ for every $m>1$ and every $t>0$.  Concerning the limit, it is proved in the above references that there exists the limit of the solutions $\{u_m(x,t): m>1\}$   as $m\to\infty$, and this limit is a function $u_\infty(x)$ that does not depend on time. If the initial function satisfies the bounds $0\le u_0(x)\le 1$, then we have
$$
\lim_{m\to\infty}u_m(x,t)=u_0(x),
$$
so that no diffusion occurs at all, cf. \cite{Benilan-Crandall81}. The more interesting case happens when $u_0$ is larger than 1 on a nontrivial set, and in that case there still exists a unique limit
$$
u_\infty(x)=\lim_{m\to\infty}u_m(x,t) \quad \mbox{ and } \quad 0\le u_\infty(x)\le 1.
$$
 This means that the upper part of the initial datum collapses at \ $t=0+$ \ to the level $u=1$, in response to the fact that the diffusivity $mu^{m-1}\to \infty$ whenever $u>1$. In other words, we are facing a {\sl singular perturbation limit} and there is an associated {\sl initial discontinuity layer}. Describing such phenomenon is the content of the mathematical theory of the mesa problem. A brief description is as follows: the upper level set  $\Omega=\{x: u_\infty(x)=1\}$ is found by solving a certain variational inequality, while  away from $\Omega$ we have $u_\infty(x)=u_0(x)$ (no diffusion takes place there).

The name of {\sl mesa problem} for this problem comes from the typical `mesa shape'  of $u_\infty(x)$ (the shapes seen in landscapes in the West of the USA). What makes the analysis more interesting is the numerically observed fact  that the mesa formation is already apparent for relatively low values of $m$, say $m\approx 6$, with typically bell-shaped initial data, cf. \cite{EHKO86}.  $\Omega$ can be much larger than the set where
$u_0(x)>1$.

We study here the limit $m\to\infty$ in the case of fractional diffusion, $0<s<1$. The analysis shows some common features, as well as quite interesting novelties worth describing. We will examine in detail some of these novel aspects. First of all, we focus on the limit behavior of  the Barenblatt solutions since this family  plays a major role in the analysis of the standard mesa problem. Here we identify the limit $m\to\infty$  of these solutions, which is a nontrivial task since they are not explicit and the limit is highly singular. We observe that, contrary to the case $s=1$, the limit $u_\infty(x,t)=F_\infty(x)$ is not compactly supported but exhibits a typical  tail with power-like spatial decay at infinity. In other words, we do not get a plain mesa in the limit but a {\sl mesa with tails}. And we are able to identify the level set $u=1$ of the limit via the solution of an obstacle problem.  To do that we have to identify two new asymptotic functions,  $G_\infty$ and $P_\infty$;  together they allow to formulate the obstacle problem. Complete proofs are given in Sections \ref{sec.lim1} and \ref{sec.obst}.

The analysis of the limit uses heavily a pair of associated functions, namely $u_m(x,t)$ and $w_m(x,t)=mu_m^m(x,t)$. Both behave very differently for large values of $m$. In the analysis of $w_m$ a second surprise arises: when applied to the Barenblatt solutions, the limit of $w_m$ as $m\to \infty$ is just the spatial profile of the self-similar solution that describes the asymptotic behaviour of a quite different nonlinear fractional diffusion model, namely
\begin{equation}\label{CVeq}
u_t=\nabla\cdot(u\nabla(-\Delta)^{-s}u)\,,
\end{equation}
studied by Caffarelli and V\'azquez in \cite{CaffVaz2}, see also Biler et al.  \cite{BIK, BKM}. This seems quite unexpected but subsequent work with Stan and Teso \cite{STVprep} shows that it is part of a more general correspondence between different models of nonlinear nonlocal diffusion.

Once this analysis is done, we devote Section \ref{sec.gen} to identify the limit for a class of nonnegative and integrable initial data. In the PME, there exists a unique limit $u_\infty(x)=\lim_{m\to\infty}u_m(x,t)$, and  $0\le u_\infty(x)\le 1$. A convenient variational inequality identifies the indicator function $w_\infty(x)$ that in turn determines the limit function $u_\infty$ and the corresponding tails with fractional type decay as $|x|\to\infty$. Our present analysis of the fractional case is only partial and will be completed  in a separate publication, but the results we present here show another difference with the standard PME: it is false that the limit functions $u_\infty(x)$ equal $u_0(x)$ at the points where $u_\infty(x)<1$.

As a final contribution of the analysis of the limit case, in Section \ref{sec.symm} we  obtain a contradiction with the standard statement of the symmetrization result (concentration comparison) that is known to be true for the standard porous medium equation, cf. \cite{Vsym82}, \cite{VANS05}. This is another remarkable difference between standard and fractional diffusion; this failure of comparison was first  demonstrated by Volzone and the author in \cite{VaVo} by completely different methods.

The analysis of the limit behaviour for general initial data is a more elaborate work that is not discussed here and we hope to perform in a future publication.

\medskip

\noindent {\sc Notations.}  $B_r(x)$ denotes the open ball in $\ren$ with center $x$ and radius $r>0$, and $\omega_N$ denotes the volume of the unit ball in $\ren$. The $s$-Laplacian operator $(-\Delta)^s$, $0<s<1$, acting in $\ren$ is precisely defined in the literature, cf. \cite{Landkof72, Stein70, Vafpme2012} among the many references. We will write $s'=1-s$. We will have to keep track of the delicate dependence of a number of constants on the values of $m \gg 1$, but we will use the same letter $C$ for different positive constants when their value is not important in the context. The dependence on $s$ will not be important in most of the text as long as $0<s<1$.

\section{Limits of Barenblatt solutions for the standard PME}\label{sec.PME}
\label{sec.fs.pme}

We re-do the analysis of the known case $s=1$ in order to introduce some detailed calculations that will fix ideas and serve as motivation. Actually, the situation for the standard porous medium equation
\begin{equation}
u_t=\Delta (u^m), \qquad m>1,
\end{equation}
posed in the whole space $x\in\ren$, $N\ge 1$, has been well-researched in the literature.  Thus, we have the following explicit formulas for the fundamental solution of the PME with data $U_m(x,0)=M\,\delta(x)$:
\begin{equation}
U_m(x,t)=t^{-\alpha}F(\xi), \qquad F_m(\xi)=(C-k\xi^2)_+^{1/(m-1)},
\end{equation}
where $\xi=x/t^{\beta}$ and
\begin{equation}
\alpha=\frac{N}{N(m-1)+2}, \quad \beta=\frac{\alpha}{N} \quad k=\frac{(m-1)\alpha}{2Nm}\,.
\end{equation}
Moreover, the free constant $C>0$ is related to the mass $M$ by the formula
$$
M=d_m\,C^{\gamma}, \qquad \gamma=\frac{N}{2(m-1)\alpha},
$$
and $d_m$ is given by the formula
$$
d_m=N\omega_N\int_0^\infty (1-ky^2)^{1/(m-1)}y^{N-1}dy=N\omega_N k^{-N/2}\int_0^\infty (1-y^2)^{1/(m-1)}y^{N-1}dy.
$$
Actually, $d_m$ depends also on the dimension $N$ but since this dependence does not play a role we will omit it as a rule. In the same we write $\alpha_m$, $\beta_m$, $\gamma_m$, and so on.

Let us now pass to the limit $m\to\infty$. We have $m \alpha_m\to 1$, $\gamma_m \to N/2$, both nontrivial limits; but $k_m\sim 1/(2Nm)\to 0$,  so we rescale $C=\widehat C/m$ and get \ $M=\widehat d_m\,\widehat C^{\gamma}$ \
with
$$
\widehat d_m=\frac{d_m}{m^\gamma}= {m^{-\gamma}}N\omega_Nk^{-N/2}\int_0^\infty (1-y^2)^{1/(m-1)}y^{N-1}dy \to\omega_N (2N)^{N/2}:=D_\infty.
$$
and put $C=\widehat C/m$  so that \ $M=\widehat d_m\,\widehat C^{\gamma}$. Using this, we easily conclude that
$$
\lim_{m\to\infty} m(U_m)^m(x,t) = W_\infty(x,tr):=\frac1{2Nt}\left(\left(\frac{M}{\omega_N}\right)^{2/N}-{|x|^2}\right)_+=\frac1{2Nt}(R_0^2-|x|^2)_+
$$
and
\begin{equation}
\lim_{m\to\infty} U_m(x,t) = U_\infty(x):=\chi_{B_{R_0}(0)}(x).
\end{equation}
We have put  $M=\omega_N R_0^N, $ and this is easily calculated on the basis that $U_\infty=1$
whenever $W_\infty>0$.

\medskip

\noindent $\bullet$
Note that $\Delta_x W_\infty(x,t)=-1/t$ \ in the set $\{(x,t): \, W_\infty>0\}$, which is the exact limit
of the well-known a priori estimate:
$$
\Delta U_m^{m-1}=-\frac{(m-1)\alpha}{mt},
$$
which holds in the same type of positivity set, $\{U_m>0\}$, for finite $m>1$.

\medskip

\noindent $\bullet$ It is interesting to write the equation for $w_m= m(u_m)^m$, which will allow us to capture part of the information in the singular limit $m\to\infty$.  The equation is
\begin{equation}
w_t=m^{1/m}w^{1-(1/m)}\Delta w\,.
\end{equation}
In the limit $w_m\to w$ it gives $w_t=w\Delta w$. This equation has $W_\infty(x,t)$ as radial
separable-variables solution, with free parameter $R_0>0$.

\medskip

\noindent $\bullet$ For the self-similar profile we have the following limits as $m\to\infty$:
\begin{equation}\label{limit.s1}
F_m(\xi)\to  F_\infty(x)=  \chi_{B_{R_0}(0)}(\xi), \qquad m(F_m(\xi))^m\to \frac1{2N}(R_0^2-|\xi|^2)_+.
\end{equation}
The limit on the left is the so-called {\sl mesa profile}. For further reference, note also that
\begin{equation}
\int_r^\infty rF(r)\,dr=(R_0-r^2)_+/2.
\end{equation}
All this is to be compared  with the calculations for the fractional case, with $0<s<1$, to be examined in the next three sections.

\section{Review of the fundamental solutions in the fractional case}\label{sec.barfrac}
\setcounter{equation}{0}

We consider next the solution $U_m(x,t)$  to the Cauchy problem \eqref{eq.ms}-\eqref{in.data} with initial data a Dirac delta, that is
$$
U_m(x,0;M) = M\delta(x), \qquad M > 0, \ m>1.
$$
This problem has been studied in \cite{vazBaren} where it is proved that for every choice of parameters $s\in (0, 1)$ and $m > m_c = \max\{(N-2s)/N, 0\}$ and every $M > 0 $ the equation admits a unique fundamental solution, which is a nonnegative continuous weak solution for $t > 0$ and takes the initial data in the sense of Radon measure, which means that
$$
\lim_{t\to 0} \int U_m(x; t)\phi(x)\, dx = M\phi(0)
$$
holds for all $\phi\in C^2_b(\ren)$. By scaling  we can reduce the study to the case $M=1$ through the formula
$$
U_m(x,t;M)=M \,U_m(x,M^{m-1}t)\,.
$$
We write in the sequel $U_m=U_m(x,t;1)$. This solution also depends on the parameters $N$ and $s$ but this dependence will be omitted as a rule since it usually plays no part in the arguments. We have the formulas
$$
U_m(x,t)=t^{-\alpha}F_m(\xi), \quad F_m(\xi) \quad \mbox{a selfsimilar profile}
$$
where \ $\xi=x/t^{\beta}$ \ and now we have the expressions
$$
\alpha=\frac{N}{N(m-1)+2s}, \quad \beta=\frac{\alpha}{N}=\frac{1}{N(m-1)+2s}\,.
$$
Moreover, $F_m$  is a bounded, positive, radial, monotone, and  H\"older continuous function that goes to zero as $|x|$ goes to  infinity.

\noindent {\sl Equation.} The self-similar profile $F=F_m$  satisfies an elliptic equation
\begin{equation}\label{sss.form}
(-\Delta )^{s}F^m =\alpha F +  \beta y\cdot \nabla F=\beta \nabla\cdot (yF)\,,
\end{equation}
so that, putting $ s'=1-s$ and integrating in $r$, we have
\begin{equation}\label{sss.form1}
\nabla (-\Delta )^{-s'} F^m=-\beta \, y\, F\,.
\end{equation}
In radial coordinates this gives
\begin{equation}\label{main.eq}
L_{s'} F^m(r)=\beta \int_r^\infty sF(s)ds\,,
\end{equation}
where $L_{s'}$ the expression of operator $(-\Delta)^{-{s'}}$ acting on radial functions. Note that the fundamental profile is a function of several parameters \ $F(r)=F_{m,s,N,M}(r)$ but only the relevant ones will be  mentioned. The scaling group acts on the profiles $F_M(r)$ for different masses $M>0$ and indeed we have
\begin{equation}\label{scaling}
F_M(r)=\mu^{2s}F_1(\mu^{1-m} r), \quad M=\mu^{N(m-1)+2s}\,,
\end{equation}
which reduces all calculations to the case $M=1$. Since $N(m-1)+2s>0$ for $M>m_c$ we get $F_M(0)\to \infty$ as $M\to\infty$. For $m\ge 1$ we also have  $\lim_{M\to\infty} F_M(r)=\infty$ for all $r>0$.

\medskip

\noindent {\sl Decay at infinity. First estimate.}
The precise behaviour of the fundamental profiles $F(y)=F_{m,s,N}(y)$ as $y\to\infty$ is a very important question in the qualitative theory. It is known in the linear case $m=1$, since $F$ is given by a linear kernel $K$ that decays like $|y|^{-(N+2s)}$, \cite{Blumenthal-Getoor}. The exact rate of decay for $m\ne 1$ is a nontrivial issue that has been carefully examined by the author in \cite{vazBaren} where it is proved  that as $r\to\infty$ we have  (al least for $m\ge 1$)
\begin{equation}
\lim_{r\to\infty} r^{N+2s}F_{m,s,N}(r)=c(m,s,N)>0\,,
\end{equation}
but this estimate in not known to be uniform in $m$ for large $m$. A less precise but uniform estimate is obtained by using the fact that $F$ is monotone as a function of $r$ and also integrable in $\ren$. Since we have the mass estimate $\int_0^\infty F_m(r)r^{N-1}dr=M/N\omega_N$ and we know that $F_m$ is monotone decreasing, we conclude that
\begin{equation}
0\le F_m(r)\le Nr^{-N}\int_0^r F_m(s)s^{N-1}ds \le \frac{M}{\omega_N r^N}.
\end{equation}
This is an upper  bound that is uniform in $m$. In the sequel we put $M=1$ without loss of generality in view of the scaling formula \eqref{scaling}.

\begin{figure}[!h]
	\begin{center}
 \includegraphics[width=0.49\textwidth,height=.45\textwidth]{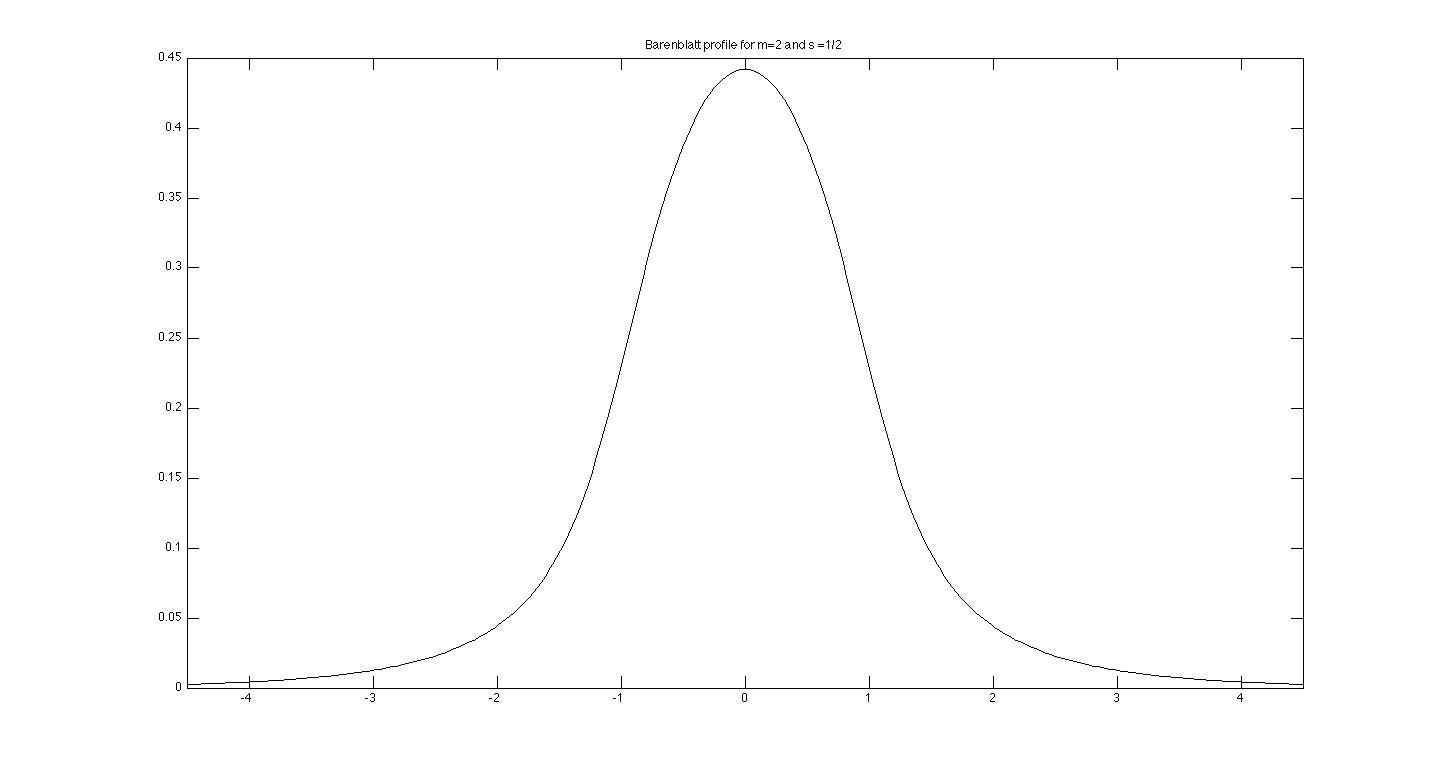}
 \includegraphics[width=0.49\textwidth,height=.45\textwidth]{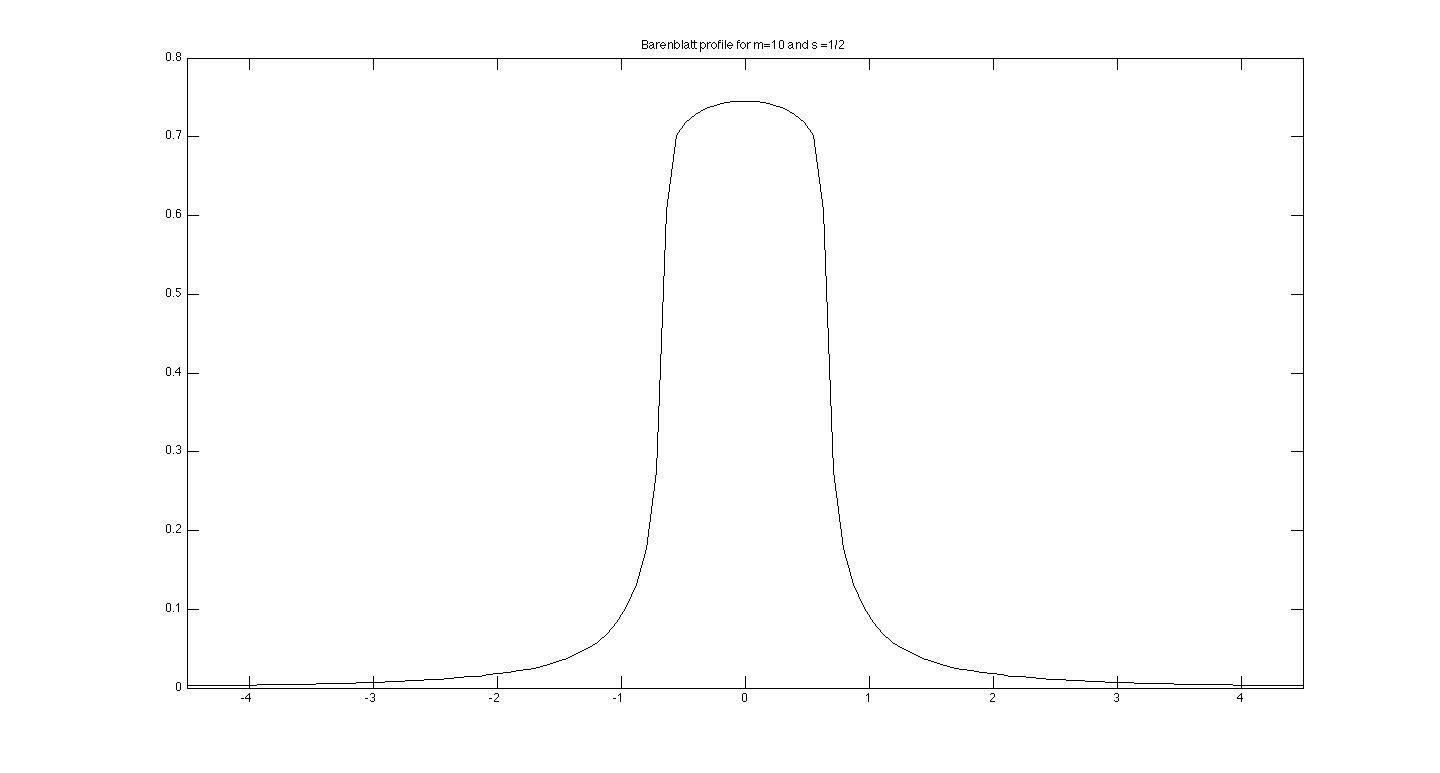}
     \end{center}
\caption{  Transition of the Barenblatt solutions to a mesa with tail in the self-similar profiles.
 Here, $N=1$, $s=0.5$ and $m=2$ (left) and $m=10$ (right).}
 \label{figure.Bar}
\end{figure}


\section{Limit of the fundamental solutions of the FPME}
\label{sec.lim1}
\setcounter{equation}{0}

We are interested in studying the limit of the family of solutions $\{U_m(r,t); m>1\}$  as $m\to\infty$. Remember that there is another parameter $s\in (0,1)$ that is kept fixed. Equivalently, we want to pass to the limit of $F_m(r)$ as $m\to\infty$ for all $r>0$. The behaviour is shown in Figure \ref{figure.Bar}.

\subsection{Limit as $m\to \infty$ for $N\ge 2$}

We begin by assuming that $N\ge 2$ since $N=1$ has some subtleties that we treat separately.
 We introduce the functions $G_m(r)=mF_m^m(r)$, that as $m\to \infty$ will have a very different behaviour compared to $F_m$. The idea of using these functions  to obtain extra information about the limit comes from the analysis Section~\ref{sec.PME}.

\begin{lem}\label{lemma4.1} Along a subsequence $m_k\to\infty$ we have $F_m(r)\to F_\infty(r)$ in $L^p(\ren)$ for all $1< p < \infty$, while  \   $G_m(r)\to G_\infty(r)$ in $L^p(\ren)$ for $1\le p< p_1=N/(N-2s)_+$. Convergence takes place also a.e. in both cases. Moreover, both limits $F_\infty$ and $G_\infty$ are non-increasing radial functions, integrable over $\ren$; $0\le F_\infty(r)\le 1$; $F_\infty(r)<1$ and $G_\infty(r)$ is zero  for $r>r_0=\omega_N^{-1/N}$;  and finally $F_\infty(r)=1$ whenever $G_\infty(r)>0$.
\end{lem}

\noindent {\sl Proof.}  (i) {\sl Passing to the limit.} Fixing some $r_1>0$ we can pass to the limit along a subsequence and we get $F_{m'}(r_1)\to F_\infty(r_1)$ as $m'=m'(r_1)\to \infty$. Doing the same for a dense countable set $\{r_k\}$ we can get the same result for all $r_k$. The limit function $F_\infty(r)$ is also nonnegative and monotone and has the same a priori bound, $\omega_N r^N\,F_\infty(r)\le {1}$. It is in principle defined for all $r=r_k$, but due to the monotonicity it can be extended to all $r\in (0,\infty)$ by limits from the left for instance, hence it is left-continuous. At all points of continuity the definition coincides with the limit $\lim_{m'\to\infty}F_m(r)$. The rest of the points (jump points) is at most countable. By Egorov's theorem $F_m\to F_\infty$ locally in all $L^p(0,\infty)$, $1\le p<\infty$, with strong convergence. At this stage it might happen that the limit $F_\infty$ will contain a Dirac delta at $r=0$,  but this will be excluded later by establishing convergence in some $L^p$ spaces, $p>1$.

\medskip

\noindent (ii) {\sl Uniform estimates  for  $G_m(r)$.} If we put $r_0=\omega_N^{-1/N}$, we have for all $r>r_0$
\begin{equation}
G_m(r)\le m\,(\omega_N r^N)^{-m}= m\,(r_0/r)^{Nm}
\end{equation}
which tends very fast to zero, uniformly on any interval of the form $[r_0+\ve,\infty)$, so that
\begin{equation}
\lim_{m\to\infty} G_m(r)=0\,, \quad r>r_0.
\end{equation}
Therefore, $G_\infty$ has compact support. Obtaining a uniform estimate on $G_m(r)$ for $r<r_0$  is more delicate and we need a different kind of argument. Taking any $r_1\in (0,r_0)$, by the monotonicity of $G_m$ w.r.t. $r$ we have  \  $ G_m(r)\ge G_m(r_1)$  for $0<r<r_1,$ so that
$$
\int_{|x|\le r_1} G_m(r)\,dx\ge \omega_N \,G_m(r_1)\,r_1^N.
$$
Using the formula for the Riesz kernel, we conclude that for all $r\ge 2r_1$ we have a constant $c_1$ (not depending on $m$) such that
$$
(-\Delta)^{-s'}G_m(x)\ge c_1\,G_m(r_1) r_1^N r^{-N+2s'}, \quad r=|x|.
$$
But on the other hand, we know the uniform bound $F_m(r)\le c_2r^{-N}$ for all large $r$, therefore
\begin{equation}\label{F.int}
\int_r^\infty sF_m(s)\,ds\le c_3r^{-N+2} \quad \mbox{for all} \ r>0.
\end{equation}
Note that for $N=2$ the last integral is just the mass in the exterior of a ball, and it is bounded above by 1. Comparing both formulas via equation \eqref{main.eq} and using the fact that $m\beta\to 1/N$ we get
$$
G_m(r_1)\, r_1^N\le c_4 r^{2s}, \quad r\ge 2r_1.
$$
In conclusion, after fixing $r=2r_1$ we get the second uniform estimate
\begin{equation}\label{Gm.sing}
G_m(r_1)\le c_5r_1^{-(N-2s)}\,,
\end{equation}
and $c_1,\dots,c_5$ do not depend on $m$.

 The two uniform estimates that we have just obtained for $G_m$ mean that $G_m(x)$ is uniformly integrable for all large $m$:
$$
\int_\ren G_m(x)\,dx\le C <\infty \qquad \mbox{for all } \ m\ge 1,
$$
and $C$ does not depend on $m$. In fact, when $m$ is large \ $G_m(x)$ is uniformly bounded in all the spaces $L^p(\ren)$ \ if \ $1\le p<p_1=N/(N-2s)$.  The details for the previous estimate  are a bit different in dimension $N=1$, see below.

\noindent (iii) {\sl Limit of \ $G_m(r)=mF_m^m(r)$ as $m\to\infty$.}   We can now apply the same argument used for $F_m$ and  prove that $G_m(r)$ converges to some $G_\infty(r)$ along some subsequences not only a.e. but also in all $L^p(\ren\setminus B_\ve(0))$, $1\le p<p_1$, with strong convergence, and
$$
G_\infty(x)\le H(|x|)=c_5|x|^{-(N-2s)}\chi_{B_{r_0}(0)}(x).
$$
This convergence eliminates the possibility of having a Dirac delta at $r=0$  in the limit. The same happens with $F_\infty(r)$.

\noindent (iv) {\sl More on the limit of  $F_m(r)$ as $m\to\infty$.} The previous results for $G_m(r)$ immediately  imply that $F_\infty(r)\le 1$ for all $r>0$.  Actually, for all $r>0$
$$
F_m^m (r)\to 0
$$
along subsequences $m'  \to \infty$. And indeed, we can say $m\to\infty$ and we do not have to take subsequences. This also means that for every $r_1>0$ there exists $m_1(r_1)$ large enough such that $F_m^m(r)< 1/2$ for $r\ge r_1$ and $m\ge m_1$, and this means that $F_m(r)\le 1$ in the same domain. On the other hand, near the origin $F_m(x)\in L^m(\ren)$ with uniform norm, hence the limit $F_m\to F_\infty$ takes place in $L^p$ for all $p< \infty$.  \qed

\medskip

Next, we establish that mass is conserved in the limit by estimating the amount of mass on the far field (what is called the tails).

\begin{lem} $F_m\to F_\infty$ in $L^1(\ren)$ and $\int F_\infty(|x|)\,dx= 1$.
 \end{lem}

\noindent {\sl Proof.}  We take a nonnegative non-increasing cutoff function $\zeta(x)$ such that $\zeta(x) = 1$ for $0 < |x|<1$, $\zeta(x) = 0 $ for $|x| > 2 $ and define $\zeta_R(x) =\zeta(x/R)$. We also put $\psi_R(x) = 1 -\zeta_R(x)$. We calculate the change in the weighted mass of the fundamental solutions $U_m$ between $t=0$ and $t=T$ for all large $m$. We take the fundamental solution with mass $M=1$ without loss of generality. We have
\begin{equation*}
\begin{array}{l}
\displaystyle 1- \int_{\ren} U_m(x,T)\zeta_R(x)\,dx=\int_{\ren} U_m(x,T)\psi_R(x)\,dx-\int_{\ren} U_m(x,0)\psi_R(x)\,dx\\[6pt]
=\displaystyle \int_0^T\int_{\ren} \partial_t U_m(x,t)\psi_R(x)\,dx dt= -\int_0^T\int_\ren (-\Delta)^s U_m^m(x,t)\psi_R(x)\,dx dt\\[6pt]
=\displaystyle -\int_0^T\int_{\ren} U_m^m(x,t)((-\Delta)^s \psi_R)(x)\,dx dt=-\int_0^T t^{N\beta-m\alpha}
\int_{\ren} F_m^m(y)((-\Delta)^s \psi_R)(yt^{\beta})\,dydt=(*)
\end{array}
\end{equation*}
We have introduced the self-similar space variable $y=xt^{-\beta}$. Of course, $U_m^m$ means $(U_m)^m$ and likewise for the notation $F_m^m$. Due to the scaling property of $\psi$
$$
((-\Delta)^s \psi_R)(y)=R^{-2s}((-\Delta)^s \psi_1)(y/R)
$$
and we also know that  $(-\Delta)^s \psi_1$ bounded in all $L^p$ spaces $1\le p\le \infty$.

\noindent (ii) We now go back to (*) to point out the estimate
$$
\begin{array}{l}
\displaystyle |\int_0^T t^{N\beta-m\alpha}
\int_\ren F_m^m(y)(-\Delta)^s \psi_R(y)\,dydt|\\[6pt]
\displaystyle  \le R^{-2s}|\int_0^T \frac{t^{N\beta-m\alpha}}{m}
\|G_m\|_1\|(-\Delta)^s \psi_1(y/R)\|_\infty\\[10pt]
\displaystyle = \frac{C\|G_m\|_1}{m(1+N\beta-m\alpha)}R^{-2s}T^{1+N\beta-m\alpha}\,.
\end{array}
$$
Recall that $\|G_m\|_1$ is uniformly bounded by Lemma \eqref{lemma4.1}. After observing that $m(1+N\beta-m\alpha)=2sm\beta\sim 2s/N$ as $m\to\infty$, we get
$$
\lim_{m\to\infty}\frac{T^{1+N\beta-m\alpha}}{m(1+N\beta-m\alpha)}\to \frac{N}{2s}
$$
so that, putting $T=1$,
$$
|\int_0^1 t^{N\beta-m\alpha}
\int_\ren F_m^m(y)(-\Delta)^s \,\psi_R(yt^{\beta})\,dydt| \le C R^{-2s}
$$
for all large $m$ and $R$, where does not depend on $R$ or $m$. Going back to the beginning of the calculation, it follows that
$$
 \left|1-\int_\ren U_m(x,1)\,\zeta_R(x)\,dx \right|\le CR^{-2s}
$$
From this we conclude the convergence of $F_m$ to $F_\infty$ in $L^1(\ren)$ and also that
$$
\int_\ren F_\infty(x)\,\zeta_R(x)\,dx=\int_\ren U_\infty(x,1)\,\zeta_R(x)\,dx=1.\mbox{\qquad \qquad \qed}
$$

\medskip

Let us now perform a further analysis of the form of $F_\infty$.

\begin{lem}   There exists $R>0$, $R<r_0$, such that $F_\infty(r)=1$ for $r<R$ and $0<F_\infty(r)<1$ for $r>R$. Moreover,
\begin{equation}
\int_0^\infty F_\infty(r)\,r^{N-1}\,dr=\frac{1}{n\omega_N}; \qquad  F_\infty(r)\sim c_1r^{-(N+2s)}\quad \mbox{as } \quad r\to \infty.
\end{equation}
\end{lem}

\noindent {\sl Proof.} (i) It follows from the previous lemma  that $F_\infty$ is monotone, and $\int F_\infty(|x|)\,dx= 1$, so that $R\le r_0=\omega_N^{-1/N}$, hence  $F_\infty$ must be less than 1 for $r>R$.

(ii) Next, we  need the equation relating the limit profiles,
\begin{equation}\label{eq.limit}
((-\Delta)^{-s'}G_\infty)(r)=\frac1{N}\int_r^\infty sF_\infty(s)ds.
\end{equation}
This is obtained by passage to the limit $m\to\infty$ in \eqref{main.eq}. The left-hand side is immediate, while for the right-hand side it comes from the Dominated Convergence Theorem if $N\ge 3$. For $N=2$ we argue as follows: the expression on the r.h.s. is just the mass of $F_\infty$ outside of the ball $B_r(0)$. Then we observe  that l.h.s gives a uniform small estimate for the mass of $F_m$ and $F_\infty(x)$ in the complement of any large ball, and we find a case of tight convergence of probability distributions.

(iii) Let us now use the equation. It is easy to prove that $(-\Delta)^{-s'}G_\infty(r)$ must be positive for all $r>0$ which means that
$\int_r^\infty s F_\infty(s)ds$ cannot be compactly supported, hence neither $F_\infty(r)$ is.
Actually, the decay rate of $(-\Delta)^{-s'}G_\infty(r)$ is $O(r^{-N-2s'})$ which means that the decay
rate of $F_\infty$ is approximately $F_\infty(r)\sim C r^{-(N+2s)}$, just as in the finite case $m< \infty$. In any case, $F_\infty(r)$ cannot be compactly supported.

(iv) We have to exclude the possibility that $R=0$ in the statement of the lemma. However,  in that case $G_\infty(r)=0$ for all $r>0$. We have to be sure that $G_\infty(r)$ is not a Dirac delta, but this has been already excluded by the convergence in some $L^p$, $p>1$. We can exclude it here in a different way: using the equation we would conclude that
$$
\int_r^\infty s F_\infty(s)ds= Cr^{-(N-2s')},
$$
which means $F_\infty(r)=c_1 r^{-(N+2s)}$ for all $r$. This contradicts the previous conclusion $F_\infty(r)\le 1$.
Summing up, $R=0$ would mean $G_\infty(r)\equiv 0$, and using \eqref{eq.limit} this would imply that $F_\infty\equiv 0$, which goes against the conservation of mass. Therefore $0<R< r_0<\infty$. \qed

\subsection{Limit of the fundamental solutions in 1D}\label{sec.lim1D}

Let us examine the proofs of this section when $N=1$. Some problems arise: thus, when $s\le 1/2$ so that $s'\ge 1/2$,  $2s'-N\ge 0$, and the argument of the Lemma  \ref{lemma4.1} has a problem at the start since the kernel involves a positive power of $|x-y|$ (or a logarithm for $s=1/2$). Moreover, even if equation \eqref{main.eq} holds for every finite $m>1$, the estimate on the asymptotic behaviour that ensures that the r.h.s. integral is finite is not uniform in $m$, and the uniform estimate we have
$\int_\RR F_m(x)dx=1$ is not sufficient.

\noindent (i) Our approach consists in taking the differentiated version \eqref{sss.form1}, i.\,e.,
\begin{equation}\label{main.eq.dim1}
-\partial_x (-\partial_{xx})^{-s'} G_m(r)=m\beta rF_m(r)\,.
\end{equation}
Using the integral kernel for $(-\partial_{xx})^{-s'}$ and differentiating we get a representation for the operator $A=-\partial_x (-\partial_{xx})^{-s'}$  (at least for smooth $f$)
$$
\begin{array}{c}
Af(x)=c\displaystyle\int_{-\infty}^x\frac{f(y)}{(x-y)^{2-2s'}}\,dy
-c\int_x^\infty\frac{f(y)}{(y-x)^{2-2s'}}\,dy\\[10pt]
=c\displaystyle\int_{-\infty}^x\frac{f(y)-f(2x-y)}{(x-y)^{2-2s'}}\,dy\,,
\end{array}$$
where $c=c(s)>0$. In this formula we have to be careful with the cancelations. By the monotonicity of $G_m$ we have nonnegative integrand for $AG_m(x)$ if $x>0$. Then,
$$
AG_m(x)\ge c\int_{-r_0}^{r_0}\frac{G_m(y)-G_m(2x-y)}{(x-y)^{2s}}\,dy=c\int_{-r_0}^{r_0}\frac{G_m(y)}{(x-y)^{2s}}\,dy
-c\int_{2x-r_0}^{2x+r_0}\frac{G_m(y)}{(x+y)^{2s}}\,dy\,.
$$
Due to the high decay rate of $G_m(x)$ for $x>2r_0$ if $m$ is very large, the last integral is very small, uniformly in $m\gg 1$ and $x>2r_0$. Hence, we conclude in the same spirit of the previous calculation for $N\ge 2$ that for all $x\ge 2r_0$ we have a constant $c_1$ (not depending on $m$) such that
$$
|A G_m(x)| \ge c_1\,(\int_{-r_0}^{r_0} G_m(x)\, dx) \, r^{-2+2s'}-\ve, \qquad r=|x|.
$$
Since on the other hand, $2xF_m(x)\le 1$ (by the integrability and monotonicity of $F_m$) we get
 the  estimate
$$
\|G_m\|_{L^1(B_{r_0})}\le c_2r^{1+2s}F_m(r)+\ve\,r^{2s}\le c_3\,r^{2s}\,, \qquad r\ge 2r_0.
$$
Fix now $r=2r_0$ to get a uniform estimate and in the limit the conclusion that $G_\infty\in L^1(B_{r_0})$, hence $G_\infty\in L^1(\RR)$.

 \noindent (ii) In order to improve that estimate we have two cases, depending on $s$ being small or not. Thus, when $1/2<s<1$ operator $A$ has symbol $-i\xi/|\xi|^{-2s'}=|\xi|^{2s-1}\mbox{\rm sign\,}(\xi)$, so the fact that $AG_m$ is bounded and that $G_m\in L^1(\RR)$ (uniformly in $m$) implies that $AG_m$ is bounded in some fractional Sobolev space and this implies that $G_m$ is uniformly in some H\"older space, and so it $G_\infty$. Note that we only need the result in a ball around the origin. In case $s=1/2$, then $A$ is a Hardy transform (but for a constant) and we conclude that $G_m$ is in $L^p$ for every $p$, uniformly in $m$, and so is $G_\infty$.

When $0<s<1/2$ we expect an estimate of the possible singularity at $x=0$ like \eqref{Gm.sing}. We argue as follows: we take a small $x>0$ and look at the kernel expression for $AG_m$ as before, but now we select the interval $x/3<y<2x/3$ to get
$$
AG_m(x)\ge c\int_{x/3}^{2x/3}\frac{G_m(y)-G_m(2x-y)}{(x-y)^{2s}}\,dy\ge c(G_m(2x/3)-G_m(4x/3))x^{1-2s}
$$
so that, using the equation $|AG_m(x)|\le c x F_m(x)\le C_1$, we get
$$
(G_m(2x/3)-G_m(4x/3))x^{1-2s}\le C_2\,.
$$
After applying this in a dyadic sequence $x_k=x_02^{-k}$ and putting $z_k=2x_k/3$ we get $G_m(z_k)\le C_3 \, z_k^{-(1-2s)}$, as we wanted to prove. We sum up the results.

\begin{lem}\label{lemma4.1.d1} The statement of Lemma {\rm \ref{lemma4.1}} is true without change for $N=1$ if $2s<1$. When $s=1/2$ there is no restriction on $p$ in the convergence of $G_m\to G_\infty$, when $s>1/2$ the convergence is uniform (and in some H\"older space). The rest of the statements holds.
\end{lem}

\begin{figure}[h!]
    \begin{center}
        \includegraphics[width=0.8\textwidth]{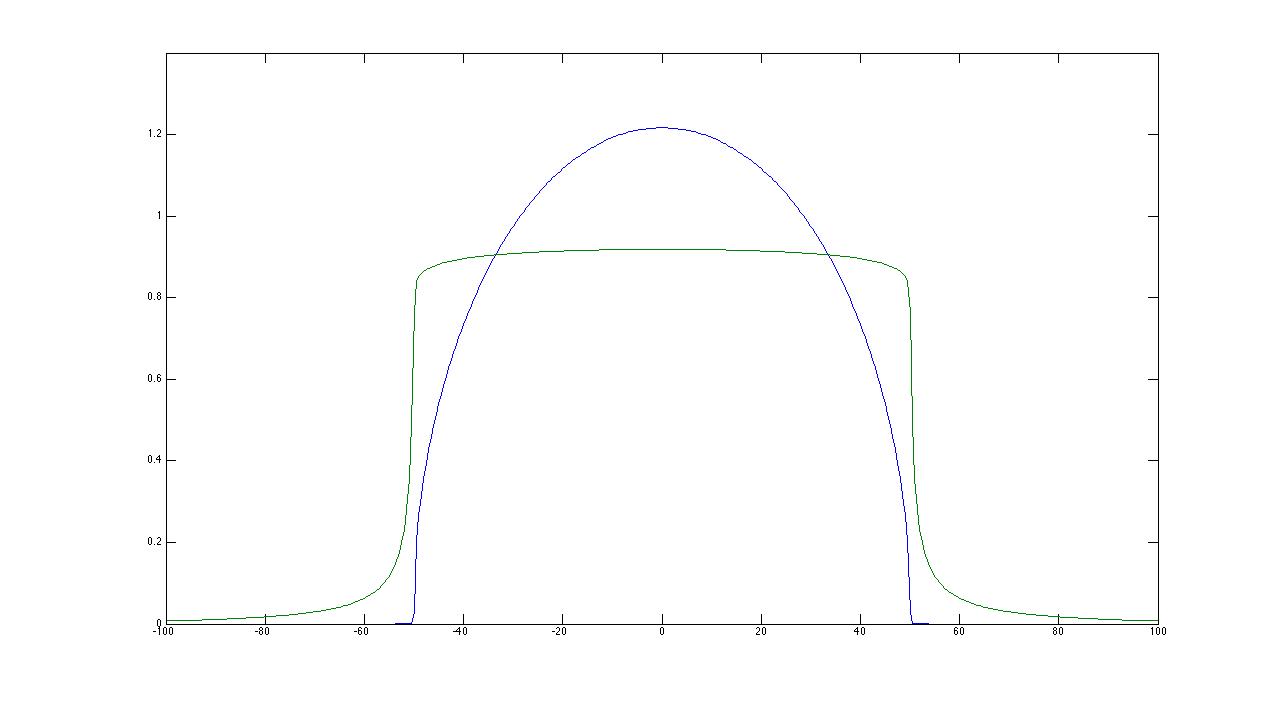}
           \end{center}
\caption{\footnotesize{Formation of the mesa shape in $F_m$ and the concave function $G_m$ for  $m=20$.
The scale on the vertical axis is real for $F_m$.}}
    \label{fig.mesaandg}
\end{figure}

Once this is established the rest of the analysis of  $G_m$ and $F_m$ of the section holds too with small changes that are not difficult. In Figure \ref{fig.mesaandg} below we represent the functions $F_m$ and $G_m$ for large  $m=20$, already showing approximation to their limit shapes.

\section{Characterization of the limit. Obstacle problem} \label{sec.obst}

The support of $G_\infty$ is some interval $[0,R_1]$ with $R_1\le R$.  We have also concluded that $G_\infty$  belongs to some $L^p$ space with $p>1$. Let us now introduce the function
\begin{equation}
P_\infty(r)=\frac1{N}\int_r^\infty sF_\infty(s)\,ds
 \end{equation}
 and let also us write without fear of confusion $P_\infty(x)=P_\infty(r)$. We know that $(-\Delta)^{s'}P_\infty(x)=G_\infty(x)$ at least in some weak sense.
Due to the form of $F_\infty$ we conclude that
\begin{equation}
P_\infty(r)=C-r^2/(2N) \quad \mbox{for } 0<r<R, \qquad P_\infty(r)\ge C- r^2/(2N) \quad \forall r>0.
\end{equation}
Hence, we know the exact shape of $P_\infty(r)$ near $r=0$ up to a constant. Moreover,  for all $r>0$ it is always above the obstacle $\Phi(r)=C-r^2/(2N)$. Here $C$ is a positive constant that may depend on the subsequence $m_k$ we have taken.

We can now present the {\bf Obstacle Problem}:  {\sl To determine radial nonnegative functions $G$ and $P$ such that \ $P(r)\ge C-r^2/(2N),$ $P(r)\to 0$ as $r\to \infty$,
\begin{equation}
 (-\Delta)^{s'}P(r)=0 \quad \mbox{\sl on the non-contact set where $P(r)> C-r^2/(2N)$}\,,
\end{equation}
 \begin{equation}
 G(r)(P(r)-\Phi(r))=0 \quad \mbox{\sl  i.\,e.,  either $G=0$ or $P$ equals the obstacle $\Phi(r):=C-(r^2/2N)$.}
 \end{equation}}
 The reader could be surprised to find  that the obstacle problem is formulated in terms of the two variables  $G$ and $P$, and not the original profile $F$. This is quite remarkable in our opinion, even if $F$ is easily obtained from $P$.

 Regarding the solution of this problem, for any given $C>0$ the Caffarelli-Silvestre theory \cite{Silv07}, \cite{CaffS}, \cite{ACS},  says that it has a solution and it is unique. The estimated regularity is $C^{1,s'}(\ren)$ for $P(x)$ and $C^{1-s'}=C^s(\ren)$ for $G(x)$.

\noindent {\bf Scaling and uniqueness.} Since our functions $P_\infty(r)$ and $G_\infty(r)$ satisfy the assumptions for some $C>0$, they coincide with the unique solution of the obstacle problem, and the stated regularity applies to them, in particular to $P_\infty$.  Then,  $rF_\infty(r)=NP_\infty'(r)\in C^{s'}$, hence  away from zero $F_\infty(r)\in C^{s'}$, while near zero $F_\infty$ is constant equal to 1.

Since the solution of the Obstacle Problem depends on the constant $C$ there is in principle an infinite family of possible solutions. The uniqueness of $C$ depends on the mass conservation law that fixes $\int F_\infty(x)\,dx=1$. Actually, when we pass from our normalized mass 1 to mass $M>0$ we easily understand what happens. We have a whole sequence of solutions of the limit problem given by the formulas
\begin{equation}\label{scaling.infty}
F_{M,\infty}(r)=F_\infty( r/M)
\end{equation}
(this is a simple scaling formula that is to be compared with \eqref{scaling}). Then, $P_{M,\infty}(r)=M^{2}P_\infty(r/M)$, so that we get all the possible constants $C_M=M^2C$, a one-to-one correspondence between mass $M$ and constant $C_M$.

We conclude from this analysis that there is a unique $C$ for which the mass of $F=F_\infty$ is one, and this ends the proof of uniqueness and implies that not only a subsequence $m'=m_k$ converges but the whole sequence $m\to\infty$ does.

Moreover, we see that the sequence of continuous and monotone decreasing functions $F_{m'}(r)$ converges to a continuous and monotone decreasing function $F_\infty(r)$, hence this convergence is  locally uniform.

\medskip

\noindent {\bf Connection with the CV fractional diffusion model. Explicit formulas.} The above obstacle problem was derived by Caffarelli and Vázquez \cite{CaffVaz2} in the study of selfsimilarity for the fractional diffusion model
\begin{equation}\label{CVeq}
u_t=\nabla\cdot(u\nabla(-\Delta)^{-\sigma/2}u)
\end{equation}
The existence and uniqueness of the self-similar solution of the form $U(x,t)=t^{-\alpha_1}F_1(xt^{-\beta_1})$
was reduced to find a solution $F_1(y)$ of that obstacle problem, and this is done via the results of  \cite{CaffS}. It is quite interesting that Biler et al. gave in \cite{BKM, BIK} an explicit formula for the solution of the evolution equation, that for the obstacle problem becomes a solution of the form
\begin{equation}
F_1(x)=(A-Bx^2)_+^{1-(\sigma/2)}
\end{equation}
with $A$ and $B$ suitable positive constants. This is based on the remarkable  explicit formula
$$
(-\Delta)^{\sigma/2}(1-|y|^2)_+^{\sigma/2}=K(\sigma,N)>0 \qquad \mbox{for} \ |y|< 1\,.
$$
due to Blumental-Getoor's \cite{Blumenthal-Getoor}  and valid for $0\le \sigma\le 2$.

In the application to our problem we must take $\sigma=2s'=2-2s$, and the solution is called $G_\infty(x)$ instead of $F_1(y)$. Since putting $R=(A/B)^{1/2}$ we get
$$
 L_{s'}(A-Bx^2)_+=AL_{s'}(1-(x/R)^2)_+=\left.AR^{-2s'}L_{s'}(1-y^2)_+\right|_{y=x/R}=-AR^{-2s'}K=- BR^{2s}K
$$
for $|x|\le R$, and since $M=cR^N$ we get the system
$$
M=cR^N, \quad BR^{2s}K=1, \qquad A=BR^2,
$$
to determine $R$, $A$ and $B$ in the explicit solution $G_\infty(x)=(A-B|x|^2)_+^{s}$. As mentioned above, the pressure is given by $P_\infty(x)=(R^2-|x|^2)/(2N)$ in the so-called coincidence set $\{P=\Phi\}$ which is the ball of radius $R$. In any case $P_\infty=(-\Delta)^{s-1}G_\infty$ in $\ren$. \qed

 \medskip

\noindent {\bf Remarks.} (1) The explicit formulas show that the positivity set of $G_\infty$ is the same as the interior of the ball where $F_\infty=1$  (i.\,e., the flat set of $F_\infty$).

\noindent  (2) It is interesting to compare the results of this section for $0<s<1$  with the explicit computations performed for the standard PME in Section \ref{sec.fs.pme}. The limit $s\to 1$ of the present results gives correct answers.  The main qualitative difference is that  $F_\infty$ is compactly supported for $s=1$, while it is {\sl not} for $s<1$. Actually,   the obstacle problem simplifies drastically when $s\to1$. Then, $s'=0$, so that $P=G$ and the alternative $G(P-\Phi)=0$ becomes $G=\Phi_+$, which is a parabola continued by zero, as we have calculated in \eqref{limit.s1}.

\noindent  (3) The difference between the initial data of the limit  process, which is $M\delta(x)$, and the value of the limit $U_m(x,t)$ for $t>0$, which is $F_\infty(x)$ is very striking, but is known in the theory of the standard mesa problems and explained as a consequence of the singular character of the limit. It takes the form of an {\sl initial discontinuity} or {\sl initial layer}. We will comment later on this issue.

\noindent  (4) The connection between the two equations that is described here has been extended to a more general correspondence in \cite{STVprep}. 


 \section{The limit for more general solutions}\label{sec.gen}

 We now consider a general initial datum given by a function $u_0\ge0$ that is bounded and integrable,  $u_0\in L^1(\ren)\cap L^\infty(\ren)$. For convenience we some times assume it  to be compactly supported too. We denote by $u_m(x,t)$ the solution of the Cauchy problem with exponent $m>1$ and fixed data $u_0$. The existence, uniqueness and properties of these solutions is studied in \cite{pqrv, pqrv2}. Figures \ref{fig.collapse} and \ref{fig.humps} below illustrate the behaviour that we expect. Note that $m=10$ produces graphs similar to $m=\infty$.

\begin{figure}[h!]
    \begin{center}
        \includegraphics[width=0.49\textwidth,height=.45\textwidth]{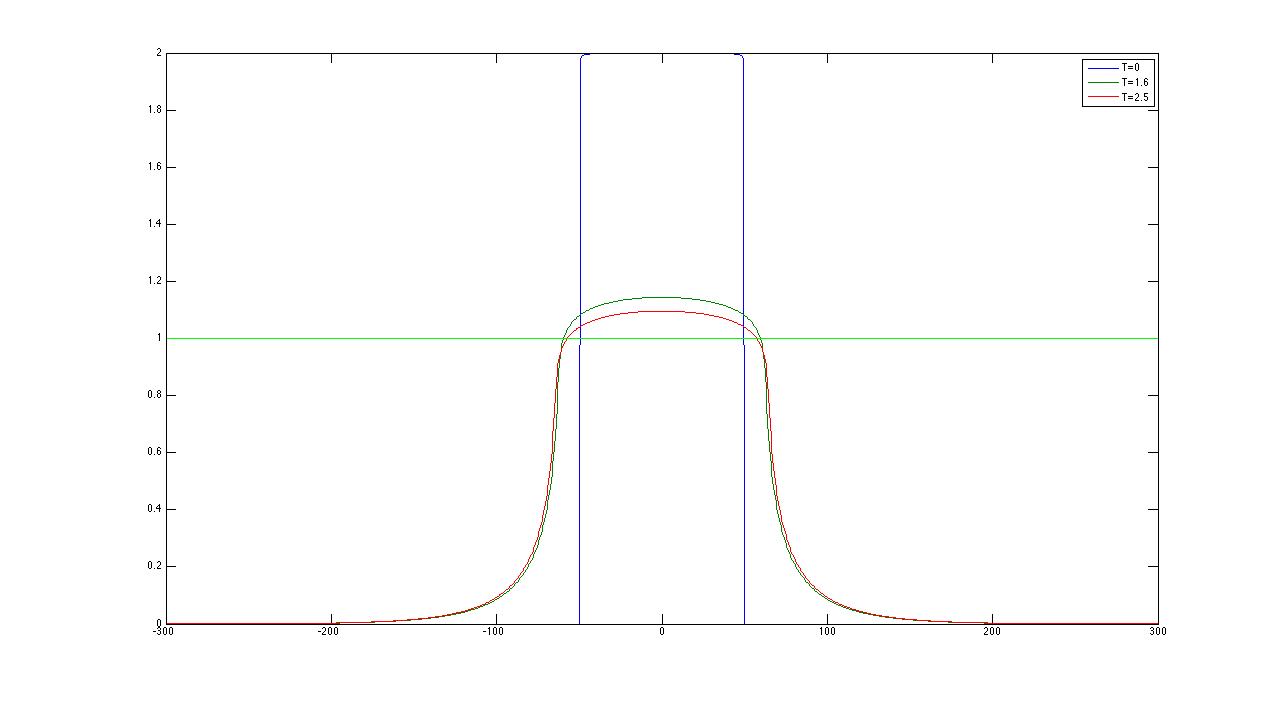}
        \includegraphics[width=0.49\textwidth,height=.45\textwidth]{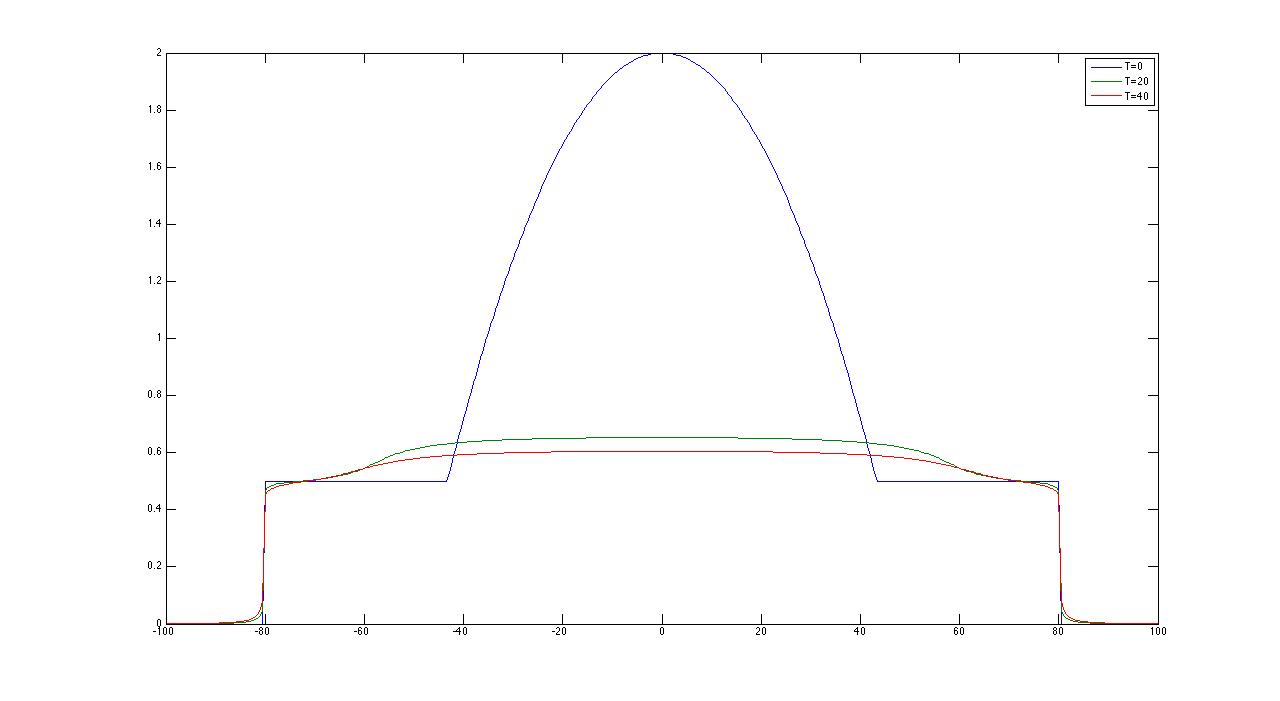}
           \end{center}
\caption{\footnotesize{Collapse of solutions to the level $u=1$. On the right, effect of a lateral step with height less than 1. In both figures $m=10$.}}
    \label{fig.collapse}
\end{figure}

\begin{figure}[h!]
    \begin{center}
        \includegraphics[width=0.49\textwidth,height=.45\textwidth]{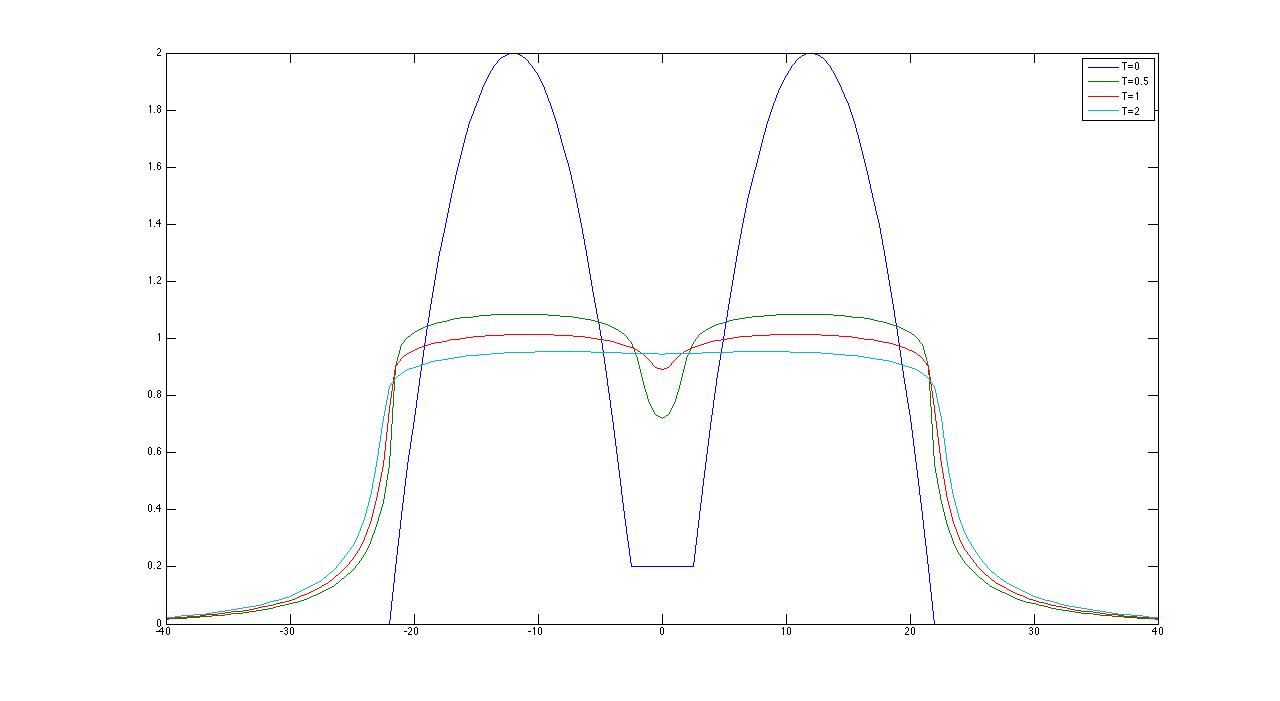}
        \includegraphics[width=0.49\textwidth,height=.45\textwidth]{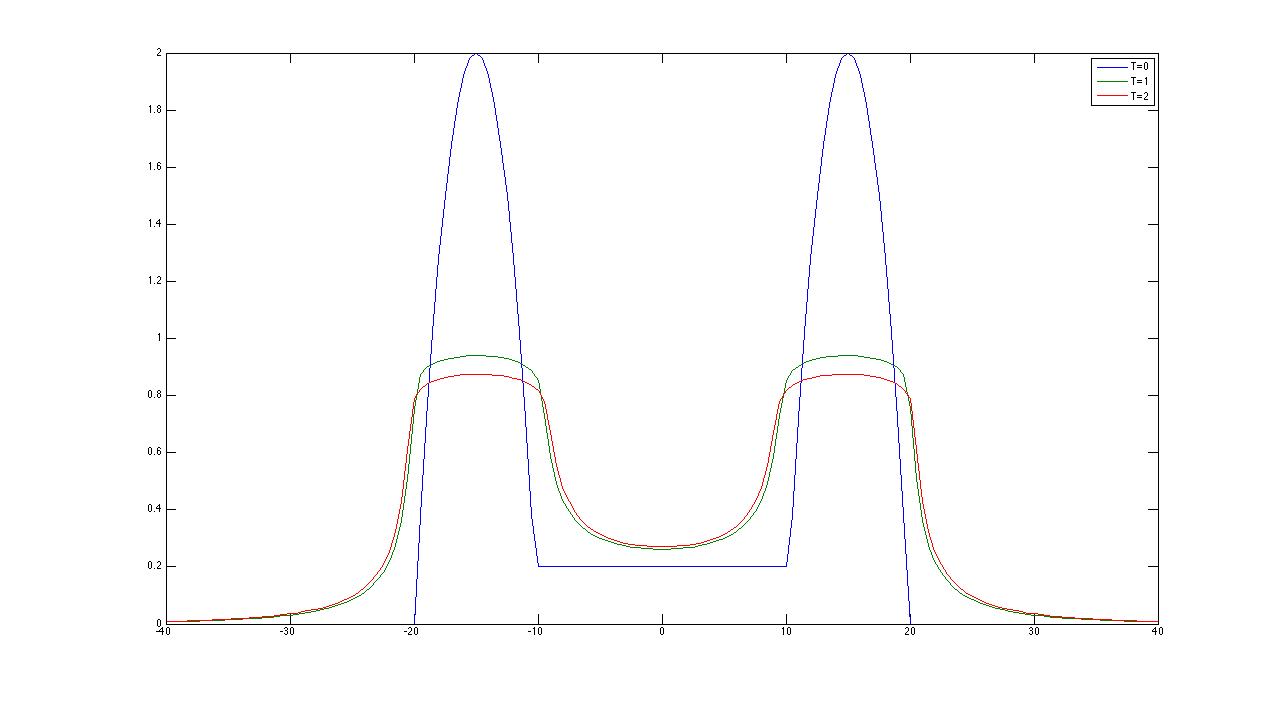}
    \end{center}
\caption{\footnotesize{Interaction of two humps at several times. Here $m=10$.}}
\label{fig.humps}
\end{figure}

\subsection{Main facts}  We want to pass to the limit in the family $\{u_m(x,t)\}_{m>1} $. The existence and properties of the limit will depend on various a priori estimates which are uniform in $m$, and will happen up to subsequences. We have to justify the type of convergence and this is what we do next.

\noindent $\bullet$ {\sl $L^p$ bounds.} First of all, it is known that for every $m>1$ and every $p\ge 1$ we have
the estimate
$$
\|u_m(\cdot,t)\|_p\le \|u_0\|_p\,.
$$
This first uniform bound allows us to pass to the limit weakly-* in $L^\infty(Q_T)$ and weakly in $L^p(\ren)$ for all $p<\infty$, along a subsequence that we denote by $m'$, to obtain a limit function $u_\infty(x,t)\in
L^\infty(Q_T)\cap L^\infty(0,T:L^1(\ren))$.

\noindent $\bullet$ {\sl Contractivity.} We have a stronger property in the $L^1(\ren)$ norm: for two solutions $u_m$ and $\widehat u_m$ with initial data $u_0$ and $\widehat u_0$
$$
    \|u_m(t)-\widehat u_m(t)\|_1\le \|u_0-\widehat u_0\|_1
$$
and this  is valid for all $m\ge 1$ and will be also valid in the limit $m\to\infty$.

\medskip

\noindent $\bullet$ {\sl Monotonicity.} It is well-known in the FPME theory \cite{BC81} that
for every $m>1$ we have for all nonnegative solutions
\begin{equation}\label{ut.1}
\partial_t u_m\ge -\frac{u_m}{(m-1)t}\,.
\end{equation}
As a consequence, we have in the limit
\begin{equation}\label{ut.1lim}
\partial_t u_\infty\ge 0\,,
\end{equation}
which means that under such general initial conditions, every limit function $u_\infty(x,t)$ is monotone nondecreasing in time.

\medskip

\noindent $\bullet$ {\sl Stronger estimate. Stationary limit.}  Moreover, for every $m>1$ we have for all nonnegative solutions with $u_0\in L^1(\ren)$ we know that
\begin{equation}\label{ut.2}
\|\partial_t u\|_1 \le \frac{2\|u_0\|_1}{(m-1)t}
\end{equation}
As a consequence, we have in the limit $\partial_t u=0$ a.e., which means that $u_\infty(x,t)$ does not depend on time for $t>0$. In other words, the limit is {\sl stationary.} This does not mean that necessarily $u_\infty(x)=u_0(x)$, because the estimate for $u_m$ is good only for $t>0$, but is singular near $t=0$. The difference between the initial data of a process and the limit of the values for the solution of the process for $t>0$ is usually labeled in the theory of singular limits as an {\sl initial discontinuity} or {\sl initial layer}. Identifying the stationary level $u_\infty(x)$ that corresponds to an initial function $u_0(x)$ is the main remaining problem of the theory.

\medskip

\noindent $\bullet$  {\sl Further regularity. Strong limit.} Estimate \eqref{ut.2} implies compactness in time for the sequence $\{u_m: m>1\}$  that we already knew to be uniformly bounded. Compactness in space depends on the $L^1$ contractivity in space.
$$
\int |u_m(x+h,t)-u_m(x,t)|\,dx\le \int |u_0(x+h)-u_0(x)|\,dx
$$
and this quantity goes to 0 as $h\to 0$. Therefore, the sequence is compact in $L^1(\ren\times (s,T))$, $0<s<T<\infty$, and the convergence $u_m\to u$ can be assumed from this moment on as being an $L^1$ convergence locally in $Q_T$, and also a convergence almost everywhere.

\medskip

\noindent $\bullet$ {\sl The case of simple limit.} There is one case in which the identification is simple, and there is no initial discontinuity.

\begin{prop} If $u_0(x)\le 1$ then $u_\infty(x,t)=u_0(x)$.
\end{prop}

\noindent {\sl Proof.} Assume first that $u_0(x)\le 1-\ve$. In that case we may write the weak solution and easily pass to the limit in the diffusion term to get
$$
\int_{\ren} u_0(x)\zeta(x)\,dx=\int_{\ren} u(x,t)\zeta(x)\,dx\,
$$
for every smooth test function $\zeta$, hence, the conclusion. For $u_0(x)\le 1$ use $L^1$ contraction. We leave these details to the reader.\qed


\subsection{New upper estimates}

At this moment we have examined two options, the last one where $u_0(x)$ is preserved in time,  and the other extreme case where $u_0(x)$ is a Dirac delta, and then there is a huge jump from this initial data to the stationary situation $u_\infty(x)$ for $t>0$. In order to examine some cases where $u_0$ is an integrable function and undergoes a jump from $t=0$ to $t=0+$ we will consider initial functions that take values larger than 1 in some nontrivial set. In order to study the limit $m\to\infty$ we need further estimates.

\noindent $\bullet$ {\sl Uniform boundedness.} This is an instance of use of the properties of the limit of the self-similar solutions that we have studied in the previous section.

\begin{prop}\label{upper.est} Suppose that $u_0(x)\ge 0$ is nonnegative, compactly supported and bounded.  Then for a.e. $x\in\ren$ and $t>0$ we have
\begin{equation}\label{gen.bound}
u_\infty(x,t)\le 1.
\end{equation}
The same is true for initial data $u_0\in L^1(\ren)$. If $u_0$ is bounded and compactly supported we have
\begin{equation}
u_m^m(x,t)\le C_1\, t^{-m\alpha}
\end{equation}
for all $t>0$ and some $C_1$ that does not depend on $m$.
\end{prop}

\noindent {\sl Proof.} (i) Let us assume that $u_0(x)\le C$ and is supported in the ball of radius $R_0$. We want to bound above the evolving solutions $u_m(x,t)$ by putting on top of them a fundamental solution with some large mass to be adjusted, and using some small shift in time. The upper bound will then be uniform in $m$ for all $m$ large enough.We consider the fundamental solution of unit mass
 $$
 U_m(x,t+\tau)=(t+\tau)^{-\alpha}F_m(x(t+\tau)^{-\beta})\,,
 $$
 and then we rescale this solution to mass $M>1$
 $$
 U_{M,m}(x,t+\tau)=M^{2s\beta}(t+\tau)^{-\alpha}F_m(x(t+\tau)^{-\beta}M^{-(m-1)\beta})\,.
 $$
We want to make sure that for some choice of $\tau=\tau_m$ and $M$ we have $U_{M,m}(R_0,0)\ge C$, i.\,e.,
 $$
 M^{2s\beta}\,\tau^{-\alpha}F_m(R_0\tau^{-\beta}M^{(m-1)\beta})\ge C\,.
 $$
Recall that for $m$ very large we have $(m-1)\beta\sim 1/N$  and $\alpha
\sim 0$.  Putting $\tau^\beta=\lambda$ and recalling that $F_m \to F_\infty$ we get sufficient conditions as follows: we first select a radius, say $R_1=R_0$ at which $F_\infty(R_1)\ge c_1$, then we put
 $$
 R_0\lambda^{-1} M^{-1/N}\le R_1 \qquad \mbox{and} \quad
 \lambda^{-N} \ge (C+\ve)/c_1.
 $$
 Therefore, select $\lambda=(2C/c_1)^{-1/N}$ and put $M=\lambda^{-N}$.
 This means  that $\tau_m= \lambda^{1/\beta}= (2C/c_1)^{-1/N\beta}\to 0$ as $m\to \infty$.

 We can now use the comparison result for the FPME. From $u_m(x,0)\le U_{M,m}(x,\tau_m)$  for all $m$ large enough  we conclude  that for all $x\in\ren$ and $t>0$
 $$
 u_m(x,t)\le U_{M,m}(x,t+\tau_m)\,,
 $$
 which in the limit gives a limit $ u_\infty(x,t)\le 1$  for every fixed $t>0$.

\medskip

(ii) It also gives immediately an upper bound on the spatial tail of the form $u_m(x,t)=O(|x|^{-(N+2s)})$, the same as in the fundamental solution.

\medskip

(iii) It is not difficult to show that the quantities $m(u_m(x,t))^m$ are uniformly bounded and integrable in $x$ by the same comparison trick. Indeed,
$$
m\|u_m(\cdot,t)\|_\infty^m \le \max_x \ mU_{M,m}(x,t+\tau_m)=mt^{-m\alpha}F_{M,m}^ m(0)=G_{M,m}(0)t^{-m\alpha}\le C_1t^{-m\alpha}.
$$
Recall now that $m\alpha\to 1$ as $m\to\infty$. Similar argument for integrability.

\medskip

(iv) It follows by approximation that estimate \eqref{gen.bound} holds for the limit solution corresponding to any initial data $u_0\ge 0$, $u_0\in L^1(\ren)$. \qed

\medskip

 \noindent {\bf Remark.} Estimate \eqref{gen.bound} is true for much more general data. We only need a bound of the form
 $$
 u_0(x)\le c+ \phi(x), \quad \mbox{with} \  c<1 \ \mbox{and} \ \phi(x)\in L^1(\ren).
 $$

We can derive a useful  consequence from the two last results and comparison.

\begin{prop} We have  $u_\infty(x,t)=1$ a.e. in the set $\{x: u_0(x)\ge 1\}$.
\end{prop}

\noindent {\sl Proof.} Define $f=\min\{u_0,1\}$. The limit of the solutions for the FPME with data $f$ is again $f$, and by comparison \ $u_\infty(x)\ge f$. Together with Proposition \ref{upper.est} it implies the result. \qed

\medskip

Therefore,  a solution with initial data $u_0$ lying somewhere above the line $u=1$ must collapse into a state $u_\infty(x)\le 1$. Since the total mass is conserved, see next, this implies that the integral of $u_\infty(x)$ on the set $\{x: u_0(x)<1\}$ must be larger than the integral of $u_0(x)$ over the same set, hence $u_\infty(x)> u_0(x)$ in a set of nonzero measure.

\medskip

\noindent $\bullet$ Next, we prove the property of mass conservation.

\begin{prop} For every $u_0\in L^1(\ren)$, $u_0\ge 0$, and every limit $u_\infty(x)$ we have
$\int u_0(x)\,dx=\int u_\infty(x)\,dx$.
\end{prop}

\noindent {\sl Proof.} We assume first that $u_0\in L^\infty(\ren)$. Using a typical cutoff function $\zeta$ and then rescaling it to $\zeta_R(x)=\zeta(Rx)$ we have
$$
\begin{array}{c}
\displaystyle |\int u_0(x)\,\zeta_R(x)\,dx-\int u_m(x,t)\,\zeta_R(x)\,dx|
=|\int_0^t\int u_m^m(-\Delta)^s\zeta_R\,dxdt|\\[8pt]
\displaystyle \le \int_0^t
 \|u_m(t)\|_\infty^{(m-1)}dt\int u_m(x,t)|(-\Delta)^s\zeta_R|\,dx=\frac{C}{2s m\beta}R^{-2s}t^{2s\beta}\,,
\end{array}
$$
where we have used the uniform bound $ \|u_m(t)\|_\infty^{(m-1)}\le Ct^{(m-1)N\beta}$ of Proposition \ref{upper.est}. In the limit it gives
$$
 |\int u_0(x)\,\zeta_R(x)\,dx-\int u_\infty(x)\,\zeta_R(x)\,dx| \le C_2R^{-2s}\,.
$$
Let now $R\to\infty$ to conclude the mass conservation rule. For general $u_0\in L^1(\ren)$ we use approximation and the property of $L^1$ contraction. \qed


\noindent {\bf Control of the initial layer.} We now introduce a new variable, $h_m(x,t)=\displaystyle\int_0^t u_m^m(x,s)ds$. In the limit $m\to\infty$ it will serve as an indicator of the initial collapse that the solution undergoes, and a locator of the resulting `debris mound', so to say. Integrating in time equation \eqref{eq.ms}, we have
\begin{equation}\label{eq.w}
(-\Delta)^s h_m(x,t)=u_0(x)-u_m(x,t)\,.
\end{equation}
In view of the a priori estimates we know that $h_m(x,t)$ converges to some $h_\infty(x)$ that does depend on $t$ for $t>0$ and we have in the limit
\begin{equation}
(-\Delta)^s h_\infty(x)=u_0(x)-u_\infty(x),
\end{equation}
that we call the $h_\infty$-equation. Since $h_\infty(x)$ has compact support (see addendum below) we conclude that $-(-\Delta)^s h_\infty$ behaves like $c|x|^{-(N+2s)}$, $c>0$ as $|x|\to \infty$. This is precisely the behaviour of $u_\infty$ under the assumption that $u_0(x)$ has compact support and $h_\infty$ is not identically zero. The last situation is implied by the assumption  $u_0(x)\ne u_\infty(x,t)$, and this in turn is true if and only if $u_0$ is not equal or less than 1. Here is the conclusion.

\begin{prop} Assume that $u_0$ is nonnegative, bounded, compactly supported. If moreover $u_0$ is not equal or less than 1 everywhere, then $h_\infty\not\equiv 0$ and
\begin{equation}
u_\infty(x)\sim c\,|x|^{-(N+2s)} \qquad \mbox{as } \ |x|\to\infty.
\end{equation}
On the other hand, if $u_0\le 1$, then $h_\infty\equiv 0$ and $u_\infty=u_0(x)$, which can have varied decay forms as $|x|\to\infty.$
\end{prop}

\medskip

\noindent {\bf Addendum. A useful computation.}   In the case of the fundamental solutions we have for $t=1$
 $$
H_m(x)=\int_0^1 t^{-Nm\beta}F_m^m(rt^{-\beta})dt= \int_1^\infty \beta^{-1}F_m^m(r\s)\s^{Nm-1}\s^{-1/\beta}d\s,
 $$
 where we have put $\s=t^{-\beta}$ (we use capital letter for the $h_m$ function of the fundamental solutions). Therefore, putting $\lambda=Nm-\beta^{-1}=N-2s$
 $$
 H_m(r)=(m\beta)^{-1}\int_1^\infty G_m(r\s)\s^{\lambda-1}d\s\le \frac1{m\beta\lambda}\,G_m(r).
 $$
 As $m\to\infty$ it converges  to $ H_\infty(r)=N \int_1^\infty G_\infty(r\s)\s^{\lambda-1}d\s$,
 which is easy to compute and has compact support. By comparison the same property of compact support is true for the $h_\infty$ corresponding to a bounded and compactly supported initial function $u_0$.

\section{Negative result for symmetrization}\label{sec.symm}

Symmetrization techniques are a very popular tool of obtaining a priori estimates for the solutions of different partial differential equations, notably those of elliptic and parabolic type. Symmetrization techniques appear in classical works like \cite{HardyLP, PS1951}. The application of Schwarz symmetrization to obtaining a priori estimates for elliptic problems is described by Weinberger in \cite{Wein62}. Sharp a priori estimates for the solutions can be derived by using comparison with a model symmetric problem. Pointwise comparison was firmly established in the works of Talenti \cite{Talenti76, Talenti79}. For parabolic problems pointwise comparison is replaced by so-called concentration comparison. In the case of the porous medium equation $u_t=\Delta u^m$ that result was established in \cite{Vsym82, VANS05}, and holds for all $m>1$. In order to state the result we want, the following definition is needed:

\noindent {\bf Definition.} Let $f,g\in L^{1}_{loc}(\ren)$ be two radially symmetric functions on $\ren$. We say that $f$ is less concentrated than $g$, and we write $f\prec g$ if for all $R>0$ we get
\begin{equation}
\int_{B_{R}(0)}f(x)\,dx\leq \int_{B_{R}(0)}g(x)\,dx.
\end{equation}

The partial order relationship $\prec$ is called \emph{comparison of mass concentrations}. The following result is well-known.

\begin{thm} Let $u_1, u_2$ be nonnegative, weak solutions of  the PME  $u_t=\Delta u^m,$ posed in $Q=\ren\times(0,\infty)$, with initial data $u_{01}, u_{02}\in L^1(\ren)\ge 0$. Assume that both $u_{02}$ and $u_{01}$ are radially symmetric and $u_{01}\prec u_{02}$. Then, for all $t>0$ we have
\begin{equation}
u_2(\cdot,t)\prec, u_1(\cdot,t).
\end{equation}
In particular, we have $\|u_2(\cdot,t)\|_p \le\|u_1(\cdot,t)\|_p$ for every $t>0$ and every $p\in [1,\infty]$.
\end{thm}

Recently, such concentration comparison has been extended by the author and Volzone \cite{VaVo} to the fractional Laplacian version $u_t+(-\Delta)^s u^m=0$ for all $m\le 1$, and the authors were surprised to find that the result {\bf does not hold} for $m>1$. We find here a confirmation for  such negative result for  the limit case $m=\infty$. As a simple consequence, it cannot hold for large enough $m$ due to the continuity of the limit demonstrated in Section \ref{eq.limit}.

\medskip

\noindent {\bf Counterexample. } It consists of radial functions. As a first candidate we take an initial function $u_{01}$ such that $0\le u_{01}(x)=2^N$  for all $|x|\le 1$,  and $u_{01}(x)=0$  otherwise. As a second candidate, we take $u_{02}(|x|)$ such that $u_{01}(x)= 1$  in a ball of radius $R=2$ and $u_{02}(x)=0$ otherwise, so that $\int u_{01}(x)\,dx=\int u_{02}(x)\,dx$, and $u_{02}\prec u_{01}$.

However, we know that $u_{\infty,2}(x)=u_{02}(x)$ is compactly supported, while $u_{\infty,1}(x)$ decays as $|x|\to\infty$ like $c|x|^{-(N+2s)}$. Therefore, it is impossible that $u_{\infty,2}\prec u_{\infty,1}$. The reader who does not like discontinuous functions will find it easy to adapt the argument and provide an example where $ u_{01}$ and $ u_{01}$ are continuous and compactly supported functions. \qed

\section*{Comments and open problems}

\noindent - We hope to continue the  analysis of the limit $m\to\infty$ for general data with the unique identification of the limit. This is  more elaborate work that involves the associated variational inequality problem.

\medskip

\noindent - We do not know what is the correct statement about comparison after symmetrization that will be valid for the solutions of the parabolic problem \eqref{eq.ms}-\eqref{in.data} and useful in the applications. Any input in this topic will be most welcome.

\

\noindent {\large\bf Acknowledgments}

\noindent  Author  partially supported by the Spanish project MTM2011-24696. He would like to thank the Isaac Newton Institute for Mathematical Sciences, Cambridge, where this work was completed during the program Free Boundary Problems and Related Topics. The author is grateful to F. del Teso for the numerical computations using the methods of \cite{Teso} and \cite{TesoVaz}. He also thanks F. Quir\'os and  B. Volzone for a careful reading of the document and useful suggestions.

\vskip 1cm


{\small


 \


\noindent {\sc Address:}

\medskip

\noindent {\sc Juan Luis V{\'a}zquez}\newline
Departamento de Matem\'{a}ticas, Universidad Aut\'{o}noma de Madrid, \\
28049 Madrid, Spain. \ e-mail: {\tt juanluis.vazquez@uam.es}

\vskip 1cm

\noindent 2000 \textit{Mathematics Subject Classification.}
35K65, 
35S10, 
26A33, 
76S05 

\medskip

\noindent \textit{Keywords and phrases.} Nonlinear fractional diffusion,
 fundamental solutions, singular limit, mesa profile, obstacle problem, symmetrization.

\end{document}